\documentclass[12pt]{article}
\usepackage{amssymb}
\include{epsf}
\newtheorem{theorem}{Theorem}
\newtheorem{corollary}{Corollary}
\def\qed{\ifhmode\unskip\nobreak\fi\ifmmode\ifinner
\else\hskip5 pt \fi\fi\hbox{\hskip5 pt \vrule width4 pt
height6 pt  depth1.5 pt \hskip 1pt }}

\begin{document}
\title{Helicoidal graphs with prescribed mean~curvature}
\author{M. Dajczer and J. H. de Lira} 
\date{}
\maketitle

\begin{abstract} We prove an existence result for helicoidal graphs
with prescribed mean curvature in a large class of warped product
spaces which comprises space forms.
\end{abstract}

\section{Helicoidal graphs}

Motivated by a question posed to one of the authors by
J. Ripoll, we address in this note the Dirichlet problem for
helicoidal graphs with prescribed mean curvature. Although the question was
originally raised for graphs in Euclidean three dimensional space, 
we solve it here in a very general context using our result in \cite{DL}.%\vspace{1ex}

In order to made precise the notion of helicoidal graph,  we first 
explain  the geometric background we consider. Then, we present the proof of the
general existence result. Finally, we discuss its particular version in the case of
Riemannian product manifolds \mbox{$M^2\times\mathbb{R}$}, that includes the
Euclidean space.\vspace{1ex}

Let $(M^n,d\sigma^2)$ be a Riemannian manifold endowed with a Killing
vector field $S_0$. Given a positive function $\varrho\in
C^\infty(M)$ so that $S_0(\varrho)=0$, we consider the warped
product manifold
$$
\bar M^{n+1}=M^n\times_\varrho\mathbb{R}
$$
with the warped metric
\begin{equation}\label{warped} 
d\bar\sigma^2 = \varrho^2 dt^2 + d\sigma^2.
\end{equation}
We denote by $S$ the lift of $S_0$ by the projection $(u,t)\in \bar
M\mapsto u\in M$. Hence, it is easy to see that $S$ is a Killing
vector field in $(\bar M,d\bar\sigma^2)$. Thus, given  
constants $a,b\in\mathbb{R}$ with $b\neq 0$, it follows that
$$
Y = aS + bT
$$
is a Killing field in $\bar M$, where $T=\partial_t$. Identifying
$M^n$ with the immersed hypersurface $M\times\{0\}\subset \bar M$, we
denote by $\Phi\colon\,\mathbb{R}\times M\to \bar M$ the flow generated by
$Y$. 

Given a bounded domain $\Omega$ in $M^n$ with boundary $\Gamma$, the
\emph{helicoidal graph}  of a function $z$ defined on $\bar\Omega$
is the hypersurface
$$
\Sigma^n =\{\Phi(z(u),u): u\in \Omega\}.
$$

The \emph{cylinder} $K$ over $\Gamma$ is the hypersurface ruled by
the flow lines of~$Y$ through $\Gamma$, i.e.,
$$
K=\{\bar u=\Phi(s,u):s\in\mathbb{R},\,u\in\Gamma\}
$$
and $H_{K}$ stands for its mean curvature when calculated inwards.\vspace{1ex}

In the following statement $\textrm{Ric}_{\bar M}$ denotes the ambient Ricci
tensor.

\begin{theorem}\label{main}
Let $\Omega$ be a $C^{2,\alpha}$ bounded domain in $M^n$. Assume
$H_{K}\ge 0$ and ${\rm Ric}_{\bar M}|_{T\Omega} \ge
-n\inf_\Gamma H_{K}^2$. Let $H\in C^\alpha(\Omega)$ and $\varphi\in
C^{2,\alpha}(\Gamma)$ be given such that $$ |H|\le \inf_\Gamma H_K.
$$
Then, there exists a unique function $z\in C^{2,\alpha}(\bar\Omega)$
satisfying $z|_\Gamma=\varphi$ whose helicoidal graph has mean
curvature function $H$ and boundary data $\varphi$.
\end{theorem}

\section{Proof of theorem 1}

Consider the submersion map $\pi\colon\,\bar M\to M$ defined by
identifying points along flow lines of $\bar Y$, i.e.,
$$
\pi(\Phi(s,u))=u, 
$$
where $s\in\mathbb{R}$ and $u\in M^n$.
Let ${\sf v}_1, \ldots, {\sf v}_n$ be a local frame defined in
$\Omega\subset M^n$. For instance, we may take the coordinate
frame relative to a choice of local coordinates in $M^n$.
Then, let $D_1,\ldots, D_n$ be the basic vector fields $\pi$-related to ${\sf v}_1,\ldots, {\sf v}_n$. It follows that
$\pi$ is a Riemannian submersion if we consider $M^n$ endowed with the
metric defined by
$$
\langle {\sf v}_i(u), {\sf v}_j(u)\rangle = \langle D_i(\Phi(s,u)),
D_j(\Phi(s,u))\rangle_{\varrho},
$$
where $\langle\cdot, \cdot\rangle_\varrho$ denotes the metric
(\ref{warped}) in $\bar M^{n+1}$.
Now, the proof  follows directly from Theorem 1 in \cite{DL}.\hfill $\qed$

\section{Helicoidal graphs in $M^2\times\mathbb{R}$}

Now we consider the case when $M^2$ is endowed with a
rotationally invariant metric. More precisely, we  have
polar coordinates $r,\theta$ such that the metric is written as
$$
ds^2=dr^2 + \psi^2(r)\,d\theta^2
$$
for some positive smooth function $\psi$. 

We define cylindrical coordinates
$r,\theta, z$ in the Riemannian product manifold $\bar M^3
=M^2\times\mathbb{R}$. We identify $M^2$ with the slice $z=0$.
Given the Killing vector field $\partial_\theta$ in $M^2$ and  $a,b\in\mathbb{R}$ with $b\neq 0$, then
$$
Y = a\partial_\theta + b\partial_z
$$
is a Killing vector field in $\bar M^3$. In terms of the cylindrical
coordinates the flow generated by $Y$ is described by
$$
\Phi_s(r,\theta,z)= (r,\theta+as,z+bs), \quad s\in \mathbb{R}.
$$

As in the general case, we define the submersion $\pi\colon\,\bar M^3\to
M^2$ by identifying points in the same orbit through a point of
$M^2$. Hence, 
$$
\pi(r,\theta, z)=(r,\theta-\frac{a}{b}z).
$$
Notice that the vector fields $\partial_r$ and
$-b\partial_\theta+a\psi^2\partial_z$ span the horizontal subspace
of $\pi$ with respect to the metric in $\bar M^3$.
Moreover, 
\begin{eqnarray}
\pi_*\partial_r = \partial_r, \quad \pi_*\partial_\theta =
\partial_\theta, \quad \pi_* \partial_z =
-\frac{a}{b}\partial_\theta,
\end{eqnarray}
what implies that the horizontal vector fields
$$
D_1 = \partial_r, \quad D_2 = \frac{b}{a^2
\psi^2+b^2}(b\partial_\theta-a\psi^2\partial_z)
$$
are $\pi$-related to the vector fields $\partial_r$ and
$\partial_\theta$ in $\mathbb{R}^2$, respectively.

Therefore, the  metric in $\bar M^3$ restricted to
horizontal subspaces has components
\begin{eqnarray*}
\langle D_1,D_1\rangle = 1, \quad \langle D_1, D_2\rangle = 0, \quad
\langle D_2, D_2\rangle = \frac{b^2 \psi^2}{a^2\psi^2+b^2}.
\end{eqnarray*}
Thus, we conclude that $\pi\colon\,\bar M^3\to M^2$ is a Riemannian
submersion if we consider in $M^2$ the metric
$$
ds^2 = dr^2 + \frac{b^2 \psi^2}{a^2\psi^2+b^2}\,d\theta^2.
$$
which coincides with the Euclidean metric when $a=0$ and $\psi(r)=r$.
\vspace{1ex}

Theorem \ref{main} in the particular case of $\mathbb{R}^3$ gives the following result.

\begin{corollary}
Let $\Omega\subset\mathbb{R}^2$ be a $C^{2,\alpha}$ bounded domain with
boundary~$\Gamma$ so that  $H_K\geq 0$. Let $H\in C^\alpha(\Omega)$
and $\varphi\in C^{2,\alpha}(\Gamma)$ be given such that $$ |H|\le
\inf_\Gamma H_K.
$$
Then, there exists a unique function $z\in C^{2,\alpha}(\bar\Omega)$
satisfying $z|_\Gamma=\varphi$ whose helicoidal graph in $\mathbb{R}^3$ has mean
curvature $H$ and boundary data~$\varphi$.
\end{corollary}

It is noteworthy to observe that this corollary reduces to Serrin's
classical existence theorem \cite{se} if we take $a=0$. We also
point out that similar results may be stated for helicoidal graphs
in hyperbolic spaces and spheres with respect to linear combinations
of Killing vector fields  generating translations along a geodesic
and rotations.

{\renewcommand{\baselinestretch}{1} \hspace*{-20ex}\begin{tabbing}
\indent \= Marcos Dajczer\\
\> IMPA \\
\> Estrada Dona Castorina, 110\\
\> 22460-320 -- Rio de Janeiro -- Brazil\\
\> marcos@impa.br\\
\end{tabbing}}

\vspace*{-4ex}

{\renewcommand{\baselinestretch}{1} \hspace*{-20ex}\begin{tabbing}
\indent \= Jorge Herbert S. de Lira\\
\> UFC - Departamento de Matematica \\
\> Bloco 914 -- Campus do Pici\\
\> 60455-760 -- Fortaleza -- Ceara -- Brazil\\
\> jorge.lira@pq.cnpq.br
\end{tabbing}}

\end{document}